\documentclass[11pt,reqno]{amsart}
\usepackage{layout}
\usepackage{latexsym,amsmath,amssymb,amsfonts,amsthm,verbatim}
\usepackage[utf8]{inputenc}
\usepackage{color}

\let \ssection=\section
\renewcommand{\section}{\setcounter{equation}{0}\ssection}

\newtheorem{Th}{Theorem}[section]
\newtheorem{theorem}[Th]{Theorem}
\newtheorem{corollary}[Th]{Corollary}
\newtheorem{proposition}[Th]{Proposition}
\newtheorem{lemma}[Th]{Lemma}
	\theoremstyle{definition}

	\theoremstyle{remark}
\newtheorem{remark}[Th]{Remark}

\def\R{{\bf R}}
 \def\E{E}
\def\D{D^{1,2}}
\def\F{\mathcal F}

\def\rv{\right|}
\def\lv{\left|}
\def\lp{\left (}
\def\rp{\right )}
\def\lc{\left [}
\def\rc{\right ]}
\def\lb{\left \lbrace}
\def\rb{\right \rbrace}
\def\e{\epsilon}

\def\dt{\frac{\partial}{\partial t}}
\def\intk{\int_{t_k}^{t_{k+1}}}

\def\tk{_{t_k}}
\def\undeux{\frac{1}{2}}

\def\lxnk{\lp X^N\tk\rp}

\def\rxnr{\lp r,X^N_r\rp}
\def\unquart{\frac{1}{4}}

\def\kk{_{k+1}}
\def\tkk{_{t_{k+1}}}
\def\inttkk{\int_t^{t\kk}}
\def\ttk{\lp t-t_k\rp}
\def\lxnkk{\lp X^N\tkk\rp}

\def\sxns{\lp s,X^N_s\rp}
\def\rvd{\right |^2}

\def\kkxnkk{\lp t\kk,X^N\tkk\rp}
\def\kxnk{\lp t_k,X^N\tk\rp}

\def\Nk{^N_{t_k}}
\def\Ntk{^N_{t_k}}
\def\Ntkk{^N_{t_{k+1}}}

\def\Nkk{^N_{t_{k+1}}}

\def\unh{\frac{1}{h}}
\def\Inkt{\mathcal I^N_k(t)}
\def\Jnkt{\mathcal J^N_k(t)}
\def\txnt{\lp t,X^N_t\rp}

\def\xnt{X^N_t}
\def\lxnt{\lp X^N_t\rp}
\def\sigmad{\sigma^2}

\def\iNk{i^N_k}
\def\jNk{j^N_k}

\def\lxntk{\lp X^N_{t_k}\rp}
\def\lxntkk{\lp X^N_{t_{k+1}}\rp}
\def\ltxnt{\lp t,X^N_t\rp}
\def\xnkk{X^N_{t_{k+1}}}
\def\xntk{X^N_{t_k}}

\def\Sh{S_h}
\def\psih{\psi_{ih}}
\def\H{{L^2(0,T)}}

\begin{document}
\title{Weak error expansion of the implicit Euler scheme}
\author{Omar Aboura}
 \email{omar.aboura@malix.univ-paris1.fr}
\address{ 
  SAMM, EA 4543,
 Universit\'e Paris 1 Panth\'eon Sorbonne,
90 Rue de Tolbiac, 75634 Paris Cedex France }
\begin{abstract}
In this paper, we extend the Talay Tubaro
theorem to the implicit Euler scheme.
\end{abstract}
\maketitle

\section{Introduction}
Let $\lp\Omega,\F, P\rp$ a probability space and
$T>0$ a fixed time.
$W$ will be a Brownian motion in $\R$ 
with respect to his own filtration $\F_t$.
We will consider the following stochastic differential equation
\begin{align}\label{diffusion}
X_t=x+\int_0^tb(X_s)ds+\int_0^t\sigma(X_s)dW_s,
\end{align}
where $x\in\R$, $b$ and $\sigma$ are real functions defined on $\R$.
It is well know that, under Lipschitz conditions on $b$ and $\sigma$,
this equation admits a unique strong solution.

For various reasons, including mathematical finance or partial
differential equations, the approximation of $\E f(X_T)$ is of importance.
One way to do this is to use an Euler scheme and to study the 
speed of convergence.
There is a vast literature on this subject and one of the pioneering work is
the paper of D.~Talay and L.~Tubaro \cite{TalayTubaro}.

Let $N\in{\bf N^*}$ and $h:=T/N$. Consider $(t_k)_{0\leq k\leq N}$ the uniform
subdivision of $[0,T]$ defined by $t_k:=kh$.
In their paper \cite{TalayTubaro} the authors deal with the explicit
Euler scheme $\lp\bar X\tk\rp_{0\leq k\leq N}$ defined as:
$\bar X_{t_0}=x$
and for $k=0,\dots,N-1$,
\begin{align}
\label{eEs}
\bar X\tkk=\bar X\tk+b\lp \bar X\tk\rp h
+
\sigma\lp \bar X\tk\rp\Delta W\kk,
\end{align}
where $\Delta W\kk:=W\tkk-W\tk$.
They study the weak error $\E f\lp \bar X_{T}\rp-\E f\lp X_T\rp$.

Here, we will use the implicit Euler scheme defined as follow:
$X^N_{t_0}=x$
and for $k=0,\dots,N-1$,
\begin{align}
\label{eq-XNtk}
X^N\tkk=X^N\tk+b\lp X^N\tkk\rp h
+
\sigma\lp X^N\tk\rp\Delta W\kk.
\end{align}
Despite the fact that this implicit scheme cannot be implemented
in most cases, it has been studied in \cite{Kloeden} but, to the best of
our knowledge, its weak error expansion has not been given.
The main reason of this study is that we believe it would be a step
in order to study a weak convergence error for SPDEs.
So far in that framework only few cases have been studied in
\cite{deBouard}-\cite{Printems} for the stochastic heat or 
Schr\"odinger equation.
\subsection*{Notations}
Let $n\in {\bf N}$ and $v,w:[0,T]\times\R\rightarrow\R$
be smooth functions. We will denote by $\partial^n v(t,x)$
the $n^{th}$ derivative of $v$ with respect to the 
space variable $x$,
except for the second derivative denoted $\Delta v(t,x)$ as usual.
Moreover, by an abuse of notation, for a function $v:\R\rightarrow\R$
and $w:[0,T]\times\R\rightarrow\R$, we will write
$(vw)(t,x):=v(x)w(t,x)$.

Given $p\in{\bf N}$, $C_p$ will denote a constant that depends
on $p$, $T$ and the coefficients $b$ and $\sigma$, but does not depend on $N$.
As usual, $C_p$ may change from line to line.

For $h$ small enough, we denote by $\Sh$ the functions defined on $\R$ by
\begin{equation}
\label{def-Sh}
\Sh(x):=1/(1-hb'(x)).
\end{equation}
It is similar to the map used by Debussche in \cite{Debussche}.
\section{The main result}
Let $u$ the (classical) solution of the following pde,
called the Kolmogorov equation:
\begin{equation}
\label{eq-edp}
\begin{cases}
\dt u(t,x)
+ b(x)\partial u(t,x)
+\undeux \sigma^2(x)\Delta u(t,x) =0,\\ 
u(T,x)=f(x).
\end{cases}
\end{equation}
The properties of $u$ will be given in the next section.
Let us mention that for $b$ and $\sigma$ smooth enough,
$u$ is smooth too.
We define the function $\psi_i:[0,T]\times\R\rightarrow\R$,
where $i$ stands for implicit, as follows for a smooth enough function $u$: 
\begin{align}\label{eq-psii}
\psi_i:=
\undeux b\partial\lp b\partial u\rp
	+\unquart\sigma^2\Delta\lp b\partial u\rp
	-\undeux b^2\Delta u
	+{1\over 8}\sigma^4\partial^4u
	-\unquart b\partial\lp\sigma^2\Delta u\rp
	-{1\over 8}\sigma^2\Delta\lp\sigma^2\Delta u\rp.
\end{align}
We are now in position to state the main result of this paper.
\begin{theorem}
\label{th-main}
Let $b,\sigma,f$ be $C^{\infty}$-functions with bounded derivatives.\\
(i) The implicit Euler scheme (\ref{eq-XNtk}) is of weak order 1,
that is, there exists a constant $C$, such that for $h$ small
enough $\lv\E f\lp X^N_T\rp-\E f\lp X_T\rp\rv\leq Ch.$\\
(ii) The weak error can be expanted as
$$\E f\lp X^N_T\rp-\E f\lp X_T\rp=h\E\int_0^T\psi_i(t,X_t)dt+O(h^2).$$
\end{theorem}
We have not given the minimal hypothesis; 
indeed we want to focus on the ideas and not on the best
set of assumptions.
The proof of this theorem is quite long;
it uses intensively the Kolmogorov equation~\eqref{eq-edp},
the It\^o and Clark-Ocone formulas.
It will be proved in the next section.
We at first compare our result with that of Talay Tubaro.
In their paper \cite{TalayTubaro}, the authors introduce the 
following function
\begin{align*}
\psi_e =&\undeux b^2\Delta u
	+\undeux b\sigma^2 \partial^3u
	+\frac{1}{8}\sigma^4\partial^4u
	+\undeux\frac{\partial^2}{\partial t^2}u
	+b\dt\partial u
	+\undeux\sigma^2\dt\Delta u,
\end{align*}
and prove the following result (see \cite{TalayTubaro} page 489).
\begin{theorem}
Let $(\bar X\tk)_{k=0,\dots,N}$ denote the explicit Euler scheme
defined by~\eqref{eEs}.
Then weak error has the following expansion
$$\E f\lp\bar X_T\rp-\E f\lp X_T\rp
	=h\E\int_0^T\psi_e(t,X_t)dt+O(h^2).$$
\end{theorem}
Applying $\dt$, $b\partial$ and finally $\undeux\sigma^2\Delta$
to \eqref{eq-edp} and summing these equations we have
\begin{align*}
{\partial^2\over \partial t^2} u
+2b\partial\dt u
+\sigma^2\Delta\dt u
=-b\partial\lp b\partial u\rp
	-\undeux b\partial\lp\sigma^2\Delta u\rp
	-\undeux\sigma^2\Delta\lp b\partial u\rp
	-\unquart\sigma^2\Delta\lp\sigma^2\Delta u\rp
\end{align*}
So we can rewrite the function $\psi_e$ as
$$\psi_e
=\undeux b^2\Delta u
	+\undeux b\sigma^2\partial^3 u
	+{1\over8}\sigma^4\partial^4u
	-\undeux b\partial\lp b\partial u\rp
	-\unquart b\partial\lp\sigma^2\Delta u\rp
	-\unquart\sigma^2\Delta\lp b\partial u\rp
	-{1\over8}\sigma^2\Delta\lp\sigma^2\Delta u\rp$$
For $b=0$, we have
$\psi_e=\psi_i={1\over8}\sigma^4\partial^4u
	-{1\over8}\sigma^2\Delta\lp\sigma^2\Delta u\rp$
as expected since in
this case the explicit and the implicit Euler scheme coincide.
We can notice that $\psi_i=\psi_e-b^2\Delta u
+\undeux\sigma^2\Delta(b\partial u)+b\partial(b\partial u)
-\undeux b\sigma^2\partial^3u$.
\section{Proof Theorem \ref{th-main}}
Here is a sketch of the proof: After
proving some property of the scheme, we 
introduce a continuous interpolation of this
scheme.
Finally, after
decomposing the weak error, we study a remainder term.
\subsection{Some tools}
\begin{proposition}[Property of $u$]
\label{prop-pte-u}
Let $\lp X^{t,x}_s\rp_{s\in[t,T]}$ denote the stochastic flow, 
that is the solution of \eqref{diffusion} starting from $x$ at time $t$
and let $u(t,x)=\E f\lp X_T^{t,x}\rp$.
Then $u$ belongs to $C^{\infty,\infty}\lp[0,T]\times\R\rp$
and satisfies the Kolmogorov equation \eqref{eq-edp}.
Moreover, for any $n,p\in{\bf N}$, there exists constants $C$ and $k$
such that
$$\lv{\partial^n\over\partial t^n}\partial^pu(t,x)\rv
	\leq C\lp 1+\lv x\rv^k\rp.$$
\end{proposition}
See for example \cite{TalayTubaro} page 486 Lemma 2.

Now we recall several results from Malliavin Calculus that will be used 
in the sequel.
For a detailled introduction, we send the reader to
D.Nualart's book \cite{Nualart}.
\begin{proposition}[Clark-Ocone formula]
\label{prop-CO}
Let $t\in[0,T]$ and $F\in L^2(\F_t)\cap\D$;
then we have for all $s\in[0,t]$
$$F=\E (F|\F_s)+\int_s^t\E\lp D_rF | \F_r\rp dW_r.$$
\end{proposition}
\begin{lemma}
\label{lem-XY-in-D}
Let $F,G\in\D$.

(i) If $F$ and $DF$ are bounded, then $FG\in\D$ and $D(FG)=FDG+GDF.$

(ii) Let $f\in C^1$ with a bounded derivative; then $f(F)\in\D$ and 
$Df(F)=f'(F)DF.$

(iii) Let $(s,t)\in[0,T]^2$ such that $s<t$ and let $F\in\D\cap L^2\lp\F_s\rp$.
Then $F\lp W_t-W_s\rp\in\D$ and
$$D_r\lc F\lp W_t-W_s\rp\rc= D_rF(W_t-W_s)
+F1_{\left\lbrace s\leq r\leq t\right\rbrace}.$$

(iv) Let $\lb H_n, n\geq 1\rb$ be a sequence of random variables in $\D$
that converges to $H$ in $L^2\lp\Omega\rp$ and such that
$\sup_n\E\lp\|DH_n\|_\H^2\rp<\infty.$
Then $H$ belongs to $\D$.
\end{lemma}
For a proof of (iii), see \cite{Nualart} Lemma 1.3.4.
Now we state some technical lemmas that will be useful in the sequel.
The following discrete Gronwall lemma is classical.
\begin{lemma}[Gronwall's lemma]
\label{lem-Gronwall}
For any nonnegative sequences $(a_k)_{0\leq k\leq N}$ and
$(b_k)_{0\leq k\leq N}$ satisfying $a_{k+1}\leq (1+Ch)a_k+b_{k+1}$,
with $C>0$. Then we have $a_k\leq e^{C(T-t_k)}\lp a_0+\sum_{i=1}^kb_i\rp$.
\end{lemma}
\begin{lemma}
\label{lem-1-Ch}
Let $L>0$; then for $h^*$ small enough (more precisely $Lh^*<1$) 
there exists $\Gamma:=\frac{L}{1-Lh^*}>0$ such that for all $h\in (0,h^*)$
we have
$\frac{1}{1-Lh}< 1+\Gamma h.$
\end{lemma}
\proof
Let  $h\in(0,h^*)$; then we have $1-Lh>1-Lh^*>0$. 
Hence $ \frac{L}{1-Lh}<\frac{L}{1-Lh^*}=\Gamma$, so that 
$Lh<\Gamma h(1-Lh)$,
which yields $ 1+\Gamma h-L h-\Gamma Lh^2=(1+\Gamma h)(1-Lh)>1$.
This concludes the proof.
\endproof
\begin{lemma}[Generalization of Young's lemma]
\label{lem-Young}
For an integer $p\geq 1$ and for $\e>0$, we have
$$(a+b)^{2^p}\leq (1+\e)^{2^p-1}a^{2^p}
	+\lp1+{1\over \e}\rp^{2^p-1}b^{2^p}.$$
\end{lemma}
\proof
We use an induction argument.
The inequality is true for $p=1$, that is
$(a+b)^2\leq(1+\e)a^2+\lp 1+{1\over \e}\rp b^2$.
Now, suppose that it is true until $p$
and will prove it for $p+1$;
indeed the induction hypothesis yields
\begin{align*}
(a+b)^{2^{p+1}}\leq &\lv (1+\e)^{2^p-1} a^{2^p}
+\lp 1+{1\over\e}\rp^{2^p-1}b^{2^p}\rvd\\
\leq & (1+\e)\lv(1+\e)^{2^p-1}\rvd \lv a^{2^p}\rvd
+\lp1+{1\over\e}\rp\lv \lp 1+{1\over\e}\rp^{2^p-1}\rvd\lv b^{2^p}\rvd.
\end{align*}
This concludes the proof.
\endproof
\subsection{Property of the implicit Euler scheme}
\begin{lemma}[Existence of the scheme]
\label{lem-xntk-exist}
For small $h$, the implicit Euler scheme \eqref{eq-XNtk} is well defined.
Moreover, for all $k=0,\ldots,N$, we have $\xntk\in L^2\lp\F\tk\rp$.
\end{lemma}
We will denote by $N_0$ the smallest integer such that 
the scheme is well defined.
\proof
For $k=0$, we have $X^N_{t_0}=x\in L^2(\F\tk)$.
Suppose that
for all $j=0,\ldots,k$, $X^N_{t_j}$ is well defined and belongs to
$L^2\lp\F_{t_j}\rp$; we prove this for $j=k+1$.
We define $\xi\kk:=X^N\tk+\sigma\lp X^N\tk\rp\Delta W\kk$.
By independence of $\Delta W\kk$ and $\F\tk$ and the linear growth of $\sigma$, 
we have that $\xi\kk\in L^2\lp\Omega\rp$.
Let $F\kk:L^2\lp \Omega\rp\rightarrow L^2\lp \Omega\rp$ be defined by
\begin{align}\label{def-Fmap}
F\kk(X):=\xi\kk+b(X)h,
\end{align}
for all $X\in  L^2\lp \Omega\rp$.
Using the Lipschitz property of $b$ we have 
$\E\lv F\kk(X)-F\kk(Y)\rvd\leq \lv \|b'\|_{\infty}h\rvd\E\lv X-Y\rvd$.
So by the fixed point theorem, if $\|b'\|_{\infty}h<1$ there exist an unique
element of $L^2(\Omega)$, noted $\xnkk$,
such that $\xnkk=F\kk\lp\xnkk\rp$. The measurability  
of $X^N\tkk$ with
respect to $\F\tkk$ is obvious.
\endproof
\begin{lemma}[Malliavin derivability]
\label{lem-xntk-in-D}
Let $h>0$ small enough; then for
all $k=0,\ldots,N$, we have $X\Nk\in \D$.
Moreover, for all $t\in[t_k,t\kk)$, we have
$D_tX^N\tkk =\Sh\lxnkk\sigma\lp X^N\tk\rp$,
where $\Sh$ is defined by \eqref{def-Sh}.
\end{lemma}
\proof
It is true for $k=0$, since $X^N_0=x$.
Now suppose that for all $j=1,\ldots,k$, $X^N_{t_j}\in\D$
and prove that $X^N\tkk\in\D$.
First, we define the following sequence in $L^2\lp\Omega\rp$:
$X^N\kk(0)=0$ and for $i\geq 0$,
$X^N\kk(i+1)=F\kk\lp X^N\kk(i)\rp$ where $F\kk$ is defined by \eqref{def-Fmap}.
Using the Lipschitz property of $F\kk$,  
since $X^N\tkk$ is a fixed point of $F\kk$, we have
$$\E\lv X\Nkk-X^N\kk(i+1)\rvd
\leq \lv \|b'\|_{\infty}h\rvd\E\lv X\Nkk-X^N\kk(i)\rvd 
\leq \lp \lv \|b'\|_{\infty}h\rvd\rp^{i+1}\E\lv X\Nkk\rvd .
$$
So $X^N\kk(i)$ converge to $X^N\tkk$ in $L^2\lp \Omega\rp$
if $\|b'\|_{\infty}h<1$.
Using the induction hypothesis, the assumptions on $\sigma$ 
and Lemma~\ref{lem-XY-in-D} (ii) and (iii), we deduce that
$\xi\kk=X^N\tk+\sigma\lp X^N\tk\rp\Delta W\kk$ belongs to $\D$.
Finally, since $b$ is Lipschitz, we deduce by induction that for all 
$i\geq 0$ $X^N_k(i)\in\D$. 
Moreover we have
$DX^N_k(i+1)=D\xi\kk+ hb'(X^N_k(i))DX^N_k(i)$
and $DX_{k+1}^N(0)=0$,
so that
$$\|DX\Nk(i+1)\|_\H^2\leq 2\|D\xi\kk\|_\H^2+2h^2\|b'\|_{\infty}^2
	\|DX\Nk(i)\|_\H^2.$$
An induction argument yields for $i\geq 1$ and $2h^2\|b'\|_{\infty}^2<1$,
$$\|DX\Nk(i)\|_\H^2\leq
	2{\|D\xi\kk\|_\H^2\over 1-2h^2\|b'\|_{\infty}^2}.
$$
Finally, we have 
$\sup_i\|DX^N_k(i)\|<\infty$.
Lemma \ref{lem-XY-in-D} (iv) proves that $X^N\tkk\in\D$.\\
Finally, let $t\in [t_k,t\kk)$; applying the Malliavin derivative $D$
to \eqref{eq-XNtk} and using Lemma \ref{lem-XY-in-D} we have
$$
D_tX^N\tkk=hb'(X^N\tkk)D_tX^N\tkk+\sigma\lp X^N\tk\rp;
$$
which concludes the proof.
\endproof
The following result gives a bound of $p^{th}$ moments of the implicit scheme.
\begin{lemma}
\label{lem-xntk-bound-moment}
Fix $p\geq 1$; then for $N_0$ large enough,
there exists a constant $C(p)>0$ such that
\begin{equation}\label{ineq-b}
\sup_{N\geq N_0} \max_{k=0,\dots,N} \E\lv X\Nk\rv^p \leq C(p).
\end{equation}
\end{lemma}
\proof
Holder's inequality shows that it suffices to consider moments 
which are power of 2, that is to check
$\sup_{N\geq N_0}\max_{k=0,\dots,N}\E\lv X^N\tk\rv^{2^p}\leq C_p$,
for every integer $p\geq 1$.
Using the generalized Young Lemma~\ref{lem-Young} 
the independence between $\Delta W_{k+1}$ and $\F\tk$,
and the fact that for all $j\in{\bf N}$,
$\E\lp\Delta W_{k+1}\rp^{2j+1}=0$, we have
for $h\in(0,h^*)$ and some constant $C_p$ depending on $h^*$
\begin{align*}
\E&\lv X\Ntkk\rv^{2^p}
\leq
\lp 1+h\rp^{2^{p}-1}
\E\lv
X\Ntk
+\sigma\lxnk\Delta W\kk
\rv^{2^p}
\\
&
+
\lp 1+\unh\rp^{2^{p}-1}
h^{2^{p}}
\E\lv b\lxnkk\rv^{2^p}
\\
\leq &
(1+C_ph)\E\lv X\Ntk\rv^{2^p}
+
\lp 1+C_ph\rp
\sum_{j=1,\dots,2^{p-1}}
{2^{p}\choose 2j}\E\lp\lv X\Ntk\rv^{2^{p}-2j}
\lv\sigma\lxnk\rv^{2j}\rp
\E\lv\Delta W\kk\rv^{2j}
\\
&
+
C_ph\lp 1+\E\lv X\Ntkk\rv^{2^p}\rp
\end{align*}
Using the identity $\E\lv \Delta W_{k+1}\rv^{2j}=C(2j)h^{j}$
and the linear growth of $\sigma$ we deduce for $h<1$
\begin{align*}
\E\lv X\Ntkk\rv^{2^p}
\leq &
\lp 1+C_ph\rp\E\lv X\Ntk\rv^{2^p}
+
C_ph
+
C_ph\E\lv X\Ntkk\rv^{2^p}
\\
&
+
(1+C_ph)C_ph
\sum_{j=1,\dots,2^{p-1}}
\E\lp\lv X\Ntk\rv^{2^{p}-2j}
\lp 1+\lv X\Ntk\rv^{2j}\rp\rp.
\end{align*}
Using the inequality:
$a^{2^{p+1}-2j}\leq a^{2^{p+1}}+1$
valid for any $a>0$,
we get for some constant $C_p>0$ and $h<1$
\begin{align*}
\E\lv X\Ntkk\rv^{2^p}
\leq &
\lp 1+C_ph\rp\E\lv X\Ntk\rv^{2^p}
+
C_ph
+
C_ph\E\lv X\Ntkk\rv^{2^p}
\end{align*}
Provided that $h$ is small enough,
the Gronwall Lemma \ref{lem-Gronwall}
and Lemma \ref{lem-1-Ch} conclude the proof.
\endproof
\def\bkn{\beta^{k,N}}
\def\zkn{z^{k,N}}
\def\gkn{\gamma^{k,N}}
\def\ekn{\eta^{k,N}}
\subsection{Some martingales and related process:
$\bkn_t, \zkn_t, \gkn_t$ and $\ekn_t$}
Let $k\in\lb0,\dots, N-1\rb$ be fixed;
in the sequel, we will use the following processes defined
for $t\in[t_k,t\kk]$ 
\begin{align}
\label{def-bzgn}
\bkn_t :=& \E\lp b\lxnkk\Big|\F_t\rp,\quad
\zkn_t:=\E\lp\left. D_tb\lxnkk \right|\F_t\rp,\\
\gkn_t:=&\sigma\lxnk+\ttk \zkn_t,\quad
\ekn_t:=\sigma\lp X^N\tk\rp\E\lp D_t
	\lp\Sh b'\rp\lp X^N\tkk\rp \Big| \F_t\rp.\nonumber
\end{align}
The following lemma describes the time evolution of these processes.

\begin{lemma} \label{lem-differential}
For all $k=0,\ldots,N-1$, and for $t\in[t_k,t\kk]$,
we have the following relation
\begin{align*}
d\bkn_t=\zkn_tdW_t,\quad
d\zkn_t=\ekn_tdW_t,\quad
d\gkn_t
= \zkn_t dt +\ekn_t(t-t_k)dW_t,\\
d\ekn_t
=\lv\sigma\lxntk\rvd\E\lp D_t\lp\Sh^3 b''\rp\lxntkk|\F_t\rp dW_t.
\end{align*}
\end{lemma}
\proof
Let $k=0,\ldots,N-1$ and let $t\in[t_k,t\kk]$.
Lemmas \ref{lem-XY-in-D} (ii), \ref{lem-xntk-exist} and \ref{lem-xntk-in-D}
and the bounds of $\|b''\|_{\infty}$ imply that $b\lxntkk\in L^2\lp\F\tkk\rp
\cap\D$ and
\begin{equation}
\label{eq-33bis}
D_tb\lp X^N\tkk\rp=
b'\lp  X^N\tkk\rp D_t X^N\tkk
=
(\Sh b')\lxnkk\sigma\lp X^N\tk\rp.
\end{equation}
So the Clark-Ocone formula in Proposition \ref{prop-CO} yields
$\bkn_t=b\lxnkk-\inttkk \zkn_sdW_s$ and hence $d\bkn_t=\zkn_tdW_t$;
where $z^{k,N}_s=\E\lp D_sb\lp X^N_{t_{k+1}}\rp | \F_s\rp$;
using \eqref{eq-33bis}
\begin{align}
\label{eq-represent-z}
\zkn_s = 
\sigma \lp X^N\tk\rp
\E\lp 
(\Sh b')\lxntkk  
\Big| \F_s\rp.
\end{align}
So taking conditionnal expectation with respect to $\F_t$, we have
\eqref{eq-represent-z}.
Since $b''$ is bounded and $b'$ Lipschitz we have that for $h$ small enough,
$\Sh b'={b'\over 1-hb'}\in C^1_b.$
So we can conclude that
$(\Sh b')\lp X^N\tkk\rp \in \D$
and using the Clark-Ocone formula we deduce that for $s\in[t_k,t\kk)$,
\begin{align*}
(\Sh b')\lp X^N\tkk\rp 
=&
\E\lp 
\Sh b'\lp X^N\tkk\rp 
| \F_s\rp
+\int_s^{t\kk}\E\lp D_u
\lc\Sh b'\lp X^N\tkk\rp \rc
| \F_u\rp dW_u,
\end{align*}
and hence
$$
d\zkn_s=\sigma\lp X^N\tk\rp\E\lp D_s
\lc{b'\over 1-hb'}\lp X^N\tkk\rp\rc
\Big| \F_s\rp
dW_s
=\eta^{k,N}_sdW_s.
$$
The differential of $\gkn_t$ is a consequence of
the previous result and It\^o's formula.
Finally, since  $\Sh b'\in C^1_b$ and $(\Sh b')'=\Sh^2b''$,
Lemma \ref{lem-xntk-in-D} implies that
$D_t(\Sh b')\lxntkk=\sigma\lxntk(\Sh^3 b'')\lxntkk$
and then 
\begin{align}\label{eq-represent-eta}
\ekn_t=\lv\sigma\lxntk\rvd\E\lp(\Sh^3 b'')\lxntkk|\F_t\rp.
\end{align}
Applying once more the Clark-Ocone formula in Proposition \ref{prop-CO}, 
we deduce 
$$\E\lp\Sh^3b''\lxntkk|\F_t\rp
=(\Sh^3b'')\lxntkk
-\int_t^{t\kk}\E \lp D_s(\Sh^3b'')\lxntkk|\F_s\rp dW_s.$$
Multiplying this by $\lv\sigma\lxntk\rvd$
and using \eqref{eq-represent-eta}, we conclude the proof.
\endproof
The next lemma provides uniform moment estimates of the above processes.
\begin{lemma}
\label{lem-ito-xnt}
Let $p\in\mathbf{N}$; then for $N_0$ large enough
there exists a constant $C_p$ such that
for $N\geq N_0$,
\begin{align*}
\max_{k=0,\dots,N-1}\sup_{t_k\leq t\leq t\kk}\E\left\lbrace\lv \bkn_t\rv^p
	+\lv \zkn_t\rv^p
	+\lv \gkn_t\rv^p
	+\lv \ekn_t\rv^p\right\rbrace
\leq C_p.
\end{align*}
\end{lemma}
\proof
Using Jensen's inequality, the Lipschitz property of $b$
and Lemma \ref{lem-xntk-bound-moment} we have
\begin{align*}
\E \lv \bkn_t\rv^p
\leq \E\lv b\lp X^N\tkk\rp\rv^p 
\leq C_p \lp 1 + \E \lv X^N\tkk\rv^p \rp\leq C_p.
\end{align*}
The identity \eqref{eq-represent-z}, Jensen's inequality,
the growth property of $\sigma$ and the 
upper estimate $b'/(1-hb')\geq C$ for small $h$,
Schwarz's inequality and 
Lemma \ref{lem-xntk-bound-moment} yield
\begin{align*}
\E\lv \zkn_t\rv^p 
	\leq &\E\lv\sigma \lp X^N\tk\rp (\Sh b')\lp X^N\tkk\rp\rv^p
	\leq C_p.
\end{align*}
Using the definition of $\gkn_t$ in \eqref{def-bzgn} 
and the previous upper estimates we deduce
$$\E\lv\gkn_t\rv^p\leq C_p\E\lv\sigma\lxnk\rv^p
	+C_ph^p\E\lv \zkn_t\rv^p
	\leq C_p.$$
Finally \eqref{eq-represent-eta}, the Jensen inequality,
the growth condition on $\sigma$, the upper bounds of
$b'$ and $b''$, Lemma \ref{lem-xntk-bound-moment}
and Schwarz's inequality yield
$$\E\lv\ekn_t\rv^p
	\leq \E\lc\lv\sigma\lxnk\rv^{2p}\lv\lp\Sh^3b''\rp\lxnkk\rv^p\rc
	\leq C_p$$
This concludes the proof.
\endproof

\subsection{Continuous interpolation}
As usual we need to introduce a continuous process 
that interpolates the implicit Euler scheme \eqref{eq-XNtk}. 
With an abuse of notation, let $\lp X^N_t\rp_{t\in[0,T]}$ be
the process defined as follow:
$X^N_0=x_0$ and for $k=0,\dots,N-1$ and $t_k\leq t\leq t\kk$
\begin{align}
X^N_t
:=
X^N\tk
+
\E\lp b\lxnkk\Big|\F_t\rp(t-t_k)
+
\sigma\lp X^N\tk\rp
\lp W_t-W\tk\rp.
\label{eq-XNt}
\end{align}
This process satisfies the following
\begin{lemma}
\label{lem-dxnt}
The process $\lp X^N_t\rp_{t\in[0,T]}$ is continuous $\F_t$-adapted and
is an interpolation of the scheme \eqref{eq-XNtk}. 
Moreover, for $k\in\lb0,\dots,N-1\rb$ and $t\in[t_k,t\kk]$ we have
$dX^N_t=\bkn_t dt+\gkn_tdW_t,$
where the process $(\bkn_t)$ and $(\gkn_t)$ are defined by \eqref{def-bzgn}.
\end{lemma}
\begin{remark}
(1) If $b=0$, \eqref{eq-XNt} corresponds to the classical interpolation
given by Talay Tubaro \cite{TalayTubaro}, since the explicit
and implicit Euler scheme are the same.\\
(2) If $b$ is linear, this continuous process differs from that used 
by Debussche in \cite{Debussche}.
Indeed, the finite dimensional analog of the interpolation correponding
to the process
$dX_t=-\beta X_tdt+\sigma\lp X_t\rp dW_t$, is defined by
$$X^D_t=X^N\tk+\int\tk^t{-\beta X^N\tk\over 1+h\beta} ds
	+\int\tk^t{\sigma\lxnk\over 1+h\beta}dW_s;$$
for $t\in[t_k,t\kk]$
(see \cite{Debussche} page 96 equation (3.2)).
In this particular case, our interpolation is given by
$$X^N_t=X^N\tk+\int\tk^t\E\lp-\beta X\tkk^N|\F_s\rp ds
	+\int\tk^t\lb\sigma\lp X^N\tk\rp
	+(s-t_k)\E\lp -D_s\beta X^N\tkk\Big|\F_s\rp\rb dW_s.$$
\end{remark}
\proof[Proof of Lemma \ref{lem-dxnt}]
The fact that $(X^N_t)$ is an $(\F_t)$-adapted process which interpolates
the scheme \eqref{eq-XNtk} is a consequence of \eqref{eq-XNt}.
The continuity is a consequence of the fact that the map
$\lp t\mapsto\E\lp X|\F_t\rp\rp$ has a continuous modification.
So, applying It\^o's formula and Lemma \ref{lem-differential}, we obtain
$d(\bkn_t(t-t_k))=\ttk \zkn_tdW_t+\bkn_t dt$,
and hence
$$dX^N_t=\bkn_t dt+\lp \sigma\lxnk+\ttk \zkn_t\rp dW_t.$$
This concludes the proof.
\endproof
We next give moment estimates of the interpolation process $X^{N}_t$.
\begin{lemma}
\label{lem-bound-xnt}
Let $p\geq 1$ and $h^*>0$ be small enough. There exits a constant $C_p>0$ 
depending on $h^*$ such that 
$$\sup_{N\geq N_0}\sup_{t\in[0,T]}\E\lv\xnt\rv^p<C_p.$$
\end{lemma}
\proof
Using Lemma \ref{lem-xntk-bound-moment}, Jensen's inequality and the
independence of $W_t-W\tk$ and $X^N\tk$, we have
for $t\in[t_k,t\kk]$:
\begin{align*}
\E\lv\xnt\rv^p
	\leq & C_p\E\lv\xntk\rv^p +C_p\E\lv\E\lp b\lxntkk|\F_t\rp\rv^p
	\lv t-t_k\rv^p+C_p\E\lv\sigma\lxntk\rv^p\lv W_t-W\tk\rv^p\\
	\leq & C_p +C_p h^p\E\lv b\lxntkk\rv^p 
	+\E\lv\sigma\lxntk\rv^p\E\lv W_t-W\tk\rv^p.
\end{align*}
Using the growth condition on $b$ and $\sigma$, moments of the normal law
and Lemma \ref{lem-xntk-bound-moment}, we deduce the result.
\endproof
The following is a straightforward
consequence of Lemmas \ref{lem-ito-xnt} and \ref{lem-bound-xnt}
\begin{corollary}
\label{cor}
Let $v:[0,T]\times\R\rightarrow\R$ be a function with polynomial growth,
and let $n_1,\ldots,n_6$ non negative integers.
Then there exists a constant $C$ independent of $h^*$ such that
for $r\in[t_k,t\kk]$ and $h\in(0,h^*)$
$$\E\lp\lv\bkn_r\rv^{n_1}\lv \zkn_r\rv^{n_2}
	\lv\gkn_r\rv^{n_3}\lv\ekn_r\rv^{n_4}
	\lv r-t_k\rv^{n_5}\lv v\rxnr\rv^{n_6}\rp\leq C.$$
\end{corollary}
\subsection{Local decomposition}
Now we return to the proof of the main theorem.
Let $u$ be the solution to the Kolmogorov equation \eqref{eq-edp}.
Using \eqref{eq-edp}, we decompose the weak error into
a sum of local errors. 
Let $\delta^N_k:=\E u\lp t\kk,X^N\tkk\rp -\E u\lp t_k,X^N\tk\rp$;
we deduce
\begin{align}
\E f\lp X^N_T\rp - \E f\lp X_T\rp
	= \E u\lp T,X^N_T\rp
	- \E u\lp 0,x\rp
	= \sum_{k=0}^{N-1}\delta^N_k.
\label{eq-loc-dec}
\end{align}
We introduce, for $t_k\leq t\leq t\kk$,
\begin{align*}
\Inkt := \E\lc\lp \bkn_t-b\lp X^N_t\rp\rp \partial u\txnt\rc,
\quad
\Jnkt := \E\lc\lp \lv\gkn_t\rvd - \sigma^2\lp X^N_t\rp\rp\Delta u\txnt\rc.
\end{align*}
Since $u\in C^{1,2}$, using It\^o's formula, Lemma \ref{lem-ito-xnt} and the
Kolmogorov equation \eqref{eq-edp} at the point $\ltxnt$,
we obtain
\begin{align}
\delta^N_k
=&
\E\intk\lb 
\dt u+\bkn_t\partial u
+\undeux\lv\gkn_t\rvd\Delta u\rb\txnt dt
\label{eq-deltaNkmoins}
\\
=&\E\intk\lb\Inkt+\undeux\Jnkt\rb dt.
\label{eq-deltaNk}
\end{align}
Now for $k=0,\ldots,N-1$, we introduce the following quantities
for $s\in[t_k,t\kk]$:
\newcommand{\ds}{{\partial\over\partial s}}
\def\ikn{i^N_k}
\def\jkn{j^N_k}
\begin{align}
\ikn(s) :=& \ds\lp b\partial u\rp\sxns
	+ \bkn_s\partial\lp b\partial u\rp\sxns\label{eqi}\\
&
	+ \undeux\lv \gkn_s\rvd\Delta\lp b\partial u\rp\sxns 
	-\bkn_s\ds\partial u\sxns\nonumber\\
&
	- \lp\lv\bkn_s\rvd+\zkn_s\gkn_s\rp\Delta u\sxns
	- \undeux\bkn_s\lv\gkn_s\rvd\partial^3u\sxns, \nonumber\\
\jkn(s) :=& \lv\gkn_s\rvd\ds\Delta u\sxns
	+ \bkn_s\lv\gkn_s\rvd\partial^3u\sxns\label{eqj}\\
&
	+ \undeux\lv \gkn_s\rv^4\partial^4u\sxns
	+ 2\gkn_s\zkn_s\Delta u\sxns
\nonumber\\
&
	+ \lv s-t_k\rvd\lv \ekn_s\rvd\Delta u\sxns
	+ 2(s-t_k)\lv \gkn_s\rvd\ekn_s\partial^3u\sxns
\nonumber\\
&
	- \ds\lp \sigma^2\Delta u\rp\sxns
	- \bkn_s\partial\lp\sigma^2\Delta u\rp\sxns
	- \undeux\lv \gkn_s\rvd\Delta\lp\sigma^2\Delta u\rp\sxns.\nonumber
\end{align}
The next two lemmas explain that, up to some sign, $\mathcal I^N_k$ (resp. $\mathcal J^N_k$)
can be viewed as an antiderivative of $\ikn$ (resp.~$\jkn$).
\begin{lemma}
\label{lem-I=inti}
For all $k=0,\ldots,N-1$, we have
$\Inkt=\E\int_t^{t\kk} \ikn(s)ds$ for $t\in[t_k,t\kk]$.
\end{lemma}
\proof
If we denote by $A:=\E\lc\bkn_t\partial u \txnt\rc
- \E \lc b\lxnkk\partial u\kkxnkk\rc $ 
and by $B:=\E\lc b\lxnkk\partial u\kkxnkk\rc
- \E\lc b\lxnt\partial u \txnt\rc $
we can write $\Inkt = A + B.$
Lemma~\ref{lem-ito-xnt}
enables us to apply It\^o's formula:
Let $v:[0,T]\times\R\rightarrow\R$ be of class $C^{1,2}$;
It\^o's formula yields
\def\dr{{\partial\over\partial r}}
\begin{align}
dv\txnt
= & \lb \dt v
	+\bkn_t\partial v
	+\undeux\lv\gkn_t\rvd\Delta v
	\rb
	\txnt dt
+ \gkn_t\partial v\txnt dW_t.
\label{eq-vrxnr}
\end{align}
Using this equation with Lemma \ref{lem-differential} we have
for $v\in C^{1,2}$ 
\begin{align}\label{eq-betav}
d\lc\bkn_rv\rxnr\rc
=&\lb\bkn_r\dr v+\lv\bkn_r\rvd\partial v+\undeux\bkn_r\lv\gkn_r\rvd\Delta v
 +\zkn_r\gkn_r\partial v\rb\rxnr dr\nonumber\\
 &+\lb\bkn_r\gkn_r\partial v+\zkn_rv\rb\rxnr dW_r.
\end{align}
The function $\Delta u$ has polynomial growth;
hence corollary \ref{cor} implies that \\
$\E\int_t^{t\kk}\lb\bkn_s\gkn_s\Delta u+\zkn_s\partial u\rb\sxns dW_s=0$.
Using equation \eqref{eq-betav} with $v=\partial u$, integrating between $t$ and $t\kk$,
using the fact that $\bkn\tkk=b\lp X^N\tkk\rp$ and
taking expectation we obtain
\begin{align}
\label{eq-I-A}
A= 
-\E\int_t^{t\kk}
\lb
\bkn_s\ds\partial u
+ \lv \bkn_s\rvd\Delta u
+ \undeux \bkn_s\lv\gkn_s\rvd\partial^3u
+ \zkn_s\gkn_s\Delta u
\rb\sxns ds.
\end{align}
Similarly, Corollary \ref{cor} implies that 
$\E\int_t^{t\kk}\gkn_s\partial\lp b\partial u\rp\sxns dW_s=0$.
Using (\ref{eq-vrxnr}) with $v=b\partial u$, 
integrating between $t$ and $t\kk$ and taking expectation yields
\begin{align*}
B= &
\E \int_t^{t\kk}
\lb
\ds\lp b\partial u\rp
+\bkn_s\partial\lp b\partial u\rp
+\undeux\lv\gkn_s\rvd\Delta\lp b\partial u\rp
\rb\sxns ds.
\end{align*}
The stochastic integral is centered by Corollary \ref{cor}.
This identity combined with \eqref{eq-I-A} concludes the proof.
\endproof
\begin{lemma}
\label{lem-J=intj}
For all $k=0,\ldots,N-1$, we have
$\Jnkt = \E\int\tk^t \jkn(s)ds$ for $t\in[t_k,t\kk]$.
\end{lemma}
\proof
Using \eqref{def-bzgn} we clearly deduce that
$\Jnkt=C+D$
where
\begin{align*}
C:=&\E\lc \lv 
\sigma\lp X^N\tk\rp 
+\lp t-t_k\rp \zkn_t\rvd
\Delta u\txnt\rc
-
\E\lc\sigma^2\lp X^N\tk\rp \Delta u\kxnk\rc,\\
D:=&\E\lc\sigma^2\lp X^N\tk\rp \Delta u\kxnk\rc
-
\E\lc\sigma^2\lxnt\Delta u\txnt\rc.
\end{align*} 
We at first rewrite the term $D$:
using (\ref{eq-vrxnr}) with $v=\sigmad\Delta u$,
integrating between $t_k$ and $t$
and taking expectation, we obtain:
\begin{align*}
D = -\E\int\tk^t\lb
\dt\lp \sigmad\Delta u\rp
+\bkn_s\partial\lp\sigmad\Delta u\rp
+\undeux\lv\gkn_s\rvd\Delta \lp \sigmad\Delta u\rp\rb\sxns 
ds,
\end{align*}
since $\sigma^2\Delta u$ has polynomial growth which implies
that the stochastic integral is centered using Corollary \ref{cor}.
It\^o's formula and Lemma \ref{lem-differential} yield
for $r\in[t_k,t\kk]$
\begin{align}
\label{eq-lvgammarvd}
d\lv\gkn_r\rvd
=\lb 2\gkn_r\zkn_r+\lv\ekn_r\rvd\lv r-t_k\rvd\rb dr
+2\gkn_r\ekn_r(r-t_k)dW_r.
\end{align}
Using this equation and \eqref{eq-vrxnr}, we have for $v$ 
of class $C^{1,2}$ and $r\in[t_k,t\kk]$
\begin{align}
\label{eq-lvgammarvdv}
d\lv\gkn_r\rvd v\rxnr
=&\lb\lv\gkn_r\rvd\dr v+\bkn_r\lv\gkn_r\rvd\partial v
+\undeux\lv\gkn_r\rv^4\Delta v \right.\\
 &\left.+2\gkn_r\zkn_r v
 +\lv\ekn_r\rvd\lv r-t_k\rvd v+2\lv\gkn_r\rvd\ekn_r(r-t_k)\partial v\rb\rxnr dr\nonumber\\
 & +\lb\lp\gkn_r\rp^3\partial v+2\gkn_r\ekn_r(r-t_k)v\rb\rxnr dW_r\nonumber.
\end{align}
Using equation \eqref{eq-lvgammarvdv} with $v=\Delta u$,
integrating between $t_k$ and $t$,
using the identity $\gkn\tk=\sigma\lxnk$
and taking expectation, we deduce
\begin{align*}
C=
\E\int\tk^t &\lb
\lv\gkn_s\rvd\ds\Delta u
+\bkn_s\lv \gkn_s\rvd\partial^3 u
+\undeux\lv\gkn_s\rv^4\partial^4u
\right.\\
&\left.
+2\gkn_s\zkn_s\Delta u
+\lv s-t_k\rvd\lv\ekn_s\rvd\Delta u
+2(s-t_k)\lv\gkn_s\rvd\ekn_s\partial^3u
\rb\sxns ds.
\end{align*}
Indeed, once more Corollary \ref{cor} and the polynomial growth
of $\partial\Delta u$ and $\Delta u$ implies that the
corresponding stochastic integral is centered.
This concludes the proof.
\endproof
Plugging the results of Lemmas \ref{lem-I=inti} and \ref{lem-J=intj}
into \eqref{eq-deltaNk} we obtain
\begin{equation}
\label{eq-dec}
\E f\lp X^N_T\rp -\E f\lp X_T\rp
=\sum_{k=0}^{N-1}\E\intk\lb
\int_t^{t\kk}\iNk(s)ds
+\undeux\int\tk^t\jNk(s)ds
\rb dt.
\end{equation}
\subsubsection*{Note:}
Thanks to Corollary \ref{cor} and the assumptions growth or boundness on
the coefficients, 
all the stochastic integrals appearing in the next section, are centered.
\subsection{Upper estimate of $\Inkt$}
We next upper estimate the difference $\phi_i(s)-\phi_i(t\kk)$,
where $\phi_i$ is one of the seven terms in the right hand side of
\eqref{eqi}
\def\tkxntk{\lp t_k,X^N\tk\rp}
\def\dtdt{{\partial^2\over \partial t^2}}
\def\dtt{{\partial^2\over \partial t^2}}
\subsubsection{The term $\phi_1(s)=\ds\lp b\partial u\rp\sxns$}
Using \eqref{eq-vrxnr} with $v=\dt(b\partial u)$,
integrating from $s$ to $t\kk$ and taking expected value
we deduce
\begin{align*}\E\ds(b\partial u)\sxns
=&\E\ds(b\partial u)\kkxnkk+R_1(s),
\end{align*}
where
\def\dss{{\partial^2\over\partial s^2}}
$$R_1(s):=-\E\int_s^{t\kk}\lb\dss(b\partial u)+\bkn_r\partial\ds(b\partial u)
	+\undeux\lv\gkn_r\rvd\Delta\ds(b\partial u)\rb\rxnr dr.$$
Futhermore, Lemmas \ref{lem-bound-xnt} and \ref{lem-ito-xnt}
and the polynomial growth of the functions involved imply that
$|R_1(s)|\leq Ch$.
\subsubsection{The term $\phi_2(s)=\bkn_s\partial\lp b\partial u\rp\sxns$}
Using \eqref{eq-betav} with $v=\partial(b\partial u)$,
integrating between $s$ and $t\kk$
and taking expectation we obtain
\begin{align*}
\E\lc\bkn_s\partial\lp b\partial u\rp\sxns\rc
	=\E \lc b\lxntkk\partial\lp b\partial u\rp\kkxnkk\rc +R_2(s),
\end{align*}
where
\begin{align*}
R_2(s):= -\E\int_s^{t\kk}\Big[&\bkn_r\lbrace\ds\partial(b\partial u)
		+\bkn_r\Delta(b\partial u)+\undeux\lv \gkn_r\rvd
		\partial^3(b\partial u)\rbrace\\
		&+\gkn_r\zkn_r\Delta(b\partial u) \Big]\rxnr dr.
\end{align*}
The polynomial growth of the functions and Lemmas \ref{lem-bound-xnt}
and \ref{lem-ito-xnt} imply that $|R_2(s)|\leq Ch$.
\subsubsection{The term $\phi_3(s)=\undeux\lv \gkn_s\rvd\Delta\lp b\partial u\rp\sxns$}
Let 
\begin{align*}&R_3(s):=
  -\undeux\E\int_s^{t\kk}
\lb\lv\gkn_r\rvd\dt\Delta(b\partial u)
+\bkn_r\lv\gkn_r\rvd\partial^3(b\partial u)
+\undeux\lv\gkn_r\rv^4\partial^4(b\partial u)
\right.\\
 &\left.+2\gkn_r\zkn_r\Delta(b\partial u)
+\lv\ekn_r\rvd\lv r-t_k\rvd\Delta(b\partial u)
 +2\lv\gkn_r\rvd\ekn_r(r-t_k)\partial^3(b\partial u) \rb\rxnr dr.
\end{align*}
Using \eqref{eq-lvgammarvdv} with $v=\undeux\Delta(b\partial u)$,
integrating between $s$ and $t\kk$, and taking expectation give us
\begin{align*}
\undeux\E\lc\lv\gkn_s\rvd\Delta\lp b\partial u\rp\sxns\rc
=\undeux\E\lc\lv \gkn\tkk\rvd\Delta\lp b\partial u\rp\kkxnkk\rc+R_3(s),
\end{align*}
with $|R_3(s)|\leq Ch$.
\subsubsection{The term $\phi_4(s)=\bkn_s\ds\partial u\sxns$}
Let
\begin{align*}
R_4(s):=\E\int_s^{t\kk}\Big[&\bkn_r\dss\partial u
		+\lv\bkn_r\rv^2\partial\ds\partial u
		+\undeux\bkn_r\lv \gkn_r\rvd\Delta\ds\partial u\\
		&+\gkn_r\zkn_r\partial\ds\partial u\Big]\rxnr dr.
\end{align*}
Using (\ref{eq-betav}) for $v=\ds\partial u$
and integrating between $s$ and $t\kk$, we obtain
\begin{align*}
-\E\lc\bkn_s\dt \partial u\sxns\rc =-\E \lc b\lxntkk\dt\partial u\kkxnkk\rc+R_4(s),
\end{align*}
with $|R_4(s)|\leq Ch$.
\subsubsection{The term $\phi_5(s)=\lv\bkn_s\rv^2\Delta u\sxns$}
Using It\^o's formula and Lemma \ref{lem-differential} we have
$$d\lv\bkn_r\rvd=\lv \zkn_r\rvd dr+2\bkn_r\zkn_rdW_r.$$
Using this equation, \eqref{eq-vrxnr} and It\^o's formula we obtain
\begin{align*}d&\lc\lv\bkn_r\rvd\Delta u\rxnr\rc=
\Big[ \lv\bkn_r\rvd\dt\Delta u+\lp\bkn_r\rp^3\partial^3 u
  +\undeux\lv\bkn_r\rvd\lv\gkn_r\rvd\partial^4u\\
  &+\lv\zkn_r\rvd\Delta u+2\bkn_r\zkn_r\gkn_r\partial^3u
\Big] \rxnr dr+dM_r,\end{align*}
where $dM_r=\lb2\bkn_r\zkn_r\Delta u
+\lv\bkn_r\rvd\gkn_r\partial^3u\rb\rxnr dW_r$
and $M_t$ is a square integrable martingale. Let
\begin{align*} R_5(s):=&\E\int_s^{t\kk}
\Big[ \lv\bkn_r\rvd\dr\Delta u+\lp\bkn_r\rp^3\partial^3 u
  +\undeux\lv\bkn_r\rvd\lv\gkn_r\rvd\partial^4u\\
  &+\zkn_r\Delta u+2\bkn_r\zkn_r\gkn_r\partial^3u
\Big] \rxnr dr.\end{align*}
Integrating between $s$ and $t\kk$
and taking expectation we have
\begin{align*}
-\E\lc\lv\bkn_s\rv^2\Delta u\sxns\rc
=-\E\lc\lv\bkn\tkk\rv^2\Delta u\kkxnkk\rc+R_5(s),
\end{align*}
with $|R_5(s)|\leq Ch$.
\subsubsection{The term $\phi_6(s)=\zkn_s\gkn_s\Delta u\sxns$}
Applying It\^o's formula to the product of
$\zkn_r\gkn_r$ and \eqref{eq-vrxnr},
and Lemma \ref{lem-differential}, we obtain for $r\in[t_k,t\kk]$
\begin{align}
\label{eq-zgv}
d&\lc\zkn_r\gkn_rv\rxnr\rc
= \lb
	\zkn_r\gkn_r\dt v
	+\zkn_r\gkn_r\bkn_r\partial v
	+\undeux \zkn_r\lp\gkn_r\rp^3\Delta v
	+\lv \zkn_r\rvd v	\right.\nonumber\\&\left.
	+\lv\ekn_r\rvd(r-t_k)v
	+\zkn_r\ekn_r(r-t_k)\gkn_r\partial v
	+\lv\gkn_r\rvd\ekn_r\partial v
   \rb\rxnr dr+dM_r		
\end{align}
where
$dM_r=  \lb
	\zkn_r\lv\gkn_r\rvd\partial v
	+\zkn_r\ekn_r(r-t_k)v
	+\gkn_r\ekn_rv
  \rb\rxnr dW_r$
and $M_r$ is a square integrable martingale.
Using equation \eqref{eq-zgv} with $v=\Delta u$,
integrating between $s$ and $t\kk$
and taking expectation give us
\begin{align*}
-\E\lc \zkn_s\gkn_s\Delta u\sxns\rc
=-\E\lc z\tkk\gkn\tkk \Delta u\kkxnkk\rc+R_6(s),
\end{align*}
with $|R_6(s)|\leq Ch$.
\subsubsection{The term $\phi_7(s)=\undeux\bkn_s\lv\gkn_s\rvd\partial^3u\sxns$}
Using Lemma \ref{lem-differential} and equation \eqref{eq-lvgammarvdv},
It\^o's formula give us for $v$ of class $C^{1,2}$
\begin{align}
\label{eq-bgv}
d&\lc\bkn_r\lv\gkn_r\rvd v\rc\rxnr
=\lb	\bkn_r\lv\gkn_r\rvd\dr v
	+\lv\bkn_r\rvd\lv\gkn_r\rvd\partial v
	+\undeux\bkn_r\lv\gkn_r\rv^4\Delta v\right.\\&\left.
	+2\bkn_r\gkn_r\zkn_r v 
	+\bkn_r\lv\gkn_r\rvd\lv r-r_k\rvd v
	+2\bkn_r\lv\gkn_r\rvd\ekn_r(r-r_k)\partial v
	+\zkn_r\lp\gkn_r\rp^3\partial v\right.\nonumber\\&\left.
	+2\gkn_r\ekn_r\zkn_r(r-r_k) v
 \rb\rxnr dr	\nonumber\\&
 +\lb	\lv\gkn_r\rvd \zkn_r v
	+\bkn_r\lp\gkn_r\rp^3\partial v
	+2\bkn_r\gkn_r\ekn_r(r-r_k)v
 \rb\rxnr dW_r.	\nonumber
\end{align}
Using this equation with $v=\undeux\partial^3u$,
integrating between $s$ and $t\kk$
and taking expectation we have
\begin{align*}
-\undeux\E\lc\bkn_s\lv\gkn_s\rvd\partial^3 u\sxns\rc
=-\undeux\E\lc\bkn\tkk\lv\gkn\tkk\rvd \partial^3 u\kkxnkk \rc+R_7(s),
\end{align*}
with $|R_7(s)|\leq Ch$.
\subsection{Upper estimate of $\Jnkt$}
\def\tphi{\tilde\phi}
We upper estimate the error $\tphi_i(s)-\tphi_i(t_k)$ where
$\tphi_i$ is one of the nine terms in the right hand side of \eqref{eqj}
\subsubsection{The term $\tphi_1(s)=\lv\gkn_s\rvd\dt\Delta u\sxns$}
\def\tR{{\tilde R}}
\def\drr{{\partial^2\over\partial r^2}}
Using \eqref{eq-lvgammarvdv} with $v=\dt\Delta u$,
integrating between $t_k$ and $s$,
taking expectation
and using the fact that $\gkn\tk=\sigma\lxnk$
we have
\begin{align*}
\E\lc\lv\gkn_s\rvd\dt\Delta u\sxns\rc
=\E\lc\lv\sigma\lxntk\rvd\dt\Delta u\tkxntk\rc+\tR_1(s),
\end{align*}
with
\begin{align*}\tR_1(s)=&\E\int\tk^s\Big[
\lv\gkn_r\rvd\drr\Delta u
	+\lbrace 2\gkn_r\zkn_r+\lv\ekn_r\rvd\lv r-t_k\rvd\rbrace
	\dr\Delta u\\
&	+\lbrace\bkn_r+2\ekn_r(r-t_k)\rbrace\lv\gkn_r\rvd
	\dr\partial^3u
	+\undeux\lv\gkn_r\rv^4\dr\partial^4u
\Big] \rxnr dr.
\end{align*}
Corollary \ref{cor} implies that $|\tR_1(s)|\leq Ch$.

\subsubsection{The term $\tphi_2(s)=\bkn_s\lv\gkn_s\rvd\partial^3u\sxns$}
For an $\F_s$-measurable random variable $Z$, we have
$\E\lp Z\bkn\tk\rp=\E\lp Zb\lp X^N\tkk\rp\rp$.
Using \eqref{eq-bgv} with $v=\partial^3 u$,
integrating between $t_k$ and $s$
and taking expectation we have
\begin{align*}
\E\lc\bkn_s\lv\gkn_s\rvd\partial^3u\sxns\rc
=\E\lc b\lxntkk\lv\sigma\lxntk\rvd\partial^3u\tkxntk\rc
	+\tR_2(s),\end{align*}
where
\begin{align*}
\tR_2(s):=\E\int\tk^s&\Big\lbrace 
	\lbrace 2\gkn_r\zkn_r\bkn_r+\bkn_r\lv\ekn_r\rvd\lv r-t_k\rvd
	+2\gkn_r\ekn_r\zkn_r(r-t_k)\rbrace\partial^3u\\
	&+\bkn_r\lv\gkn_r\rvd\dr\partial^3u
	+\undeux\bkn_r\lv\gkn_r\rv^4\partial^5u\\
	&+\lbrace\lv\bkn_r\rvd+2\bkn_r\ekn_r(r-t_k)
		+\gkn_r\zkn_r
	\rbrace\lv\gkn_r\rvd\partial^4u\Big\rbrace \rxnr dr.
\end{align*}
Corollary \ref{cor} implies that
$|\tR_2(s)|\leq Ch$.

\subsubsection{The term $\tphi_3(s)=\undeux\lv\gkn_s\rv^4\partial^4u\sxns$}
Using Lemma \ref{lem-differential} and It\^o's formula 
we deduce
\begin{align*}
d\lv\gkn_r\rv^4
=\lv\gkn_r\rvd\lbrace 4\gkn_r\zkn_t+6\lv\ekn_r\rvd\lv r-t_k\rvd\rbrace dr
+4\lp\gkn_r\rp^3\ekn_r(r-t_k)dW_r.
\end{align*}
Using this equation and \eqref{eq-vrxnr} with $v=\undeux\partial^4u$
and applying It\^o formula, we have
\begin{align*}
\undeux\E\lc\lv\gkn_s\rv^4\partial^4u\sxns\rc
=\undeux\E\lc\lv\sigma\lxntk\rv^4\partial^4u\tkxntk\rc+\tR_{3}(s),
\end{align*}
where $\tR_{3}(s)\leq Ch$ by Corollary \ref{cor}.
\subsubsection{The term $\tphi_4(s)=2\gkn_s\zkn_s\Delta u\sxns$}
Using \eqref{eq-zgv} with $v=2\Delta u$ we have 
$\E\lc2\gkn_s\zkn_s\Delta u\sxns\rc=\E\lc2\gkn\tk z\tk\Delta u\kxnk\rc+\tR_{4}(s)$,
and Corollary \ref{cor} implies
$|\tR_{4}(s)|\leq Ch$.
\subsubsection{The term $\tphi_5(s):=\tR_{5}(s):=\lv s-t_k\rvd\lv \ekn_s\rvd\Delta u\sxns+
2(s-t_k)\lv \gkn_s\rvd\ekn_s\partial^3u\sxns$}
Using Corollary~\ref{cor}, we have $|\tR_{5}(s)|\leq Ch$.
\subsubsection{The term $\tphi_6(s)=\dt\lp \sigma^2\Delta u\rp\sxns$}
Using \eqref{eq-vrxnr} with $v=\dt\lp \sigma^2\Delta u\rp$,
integrating between $t_k$ and $s$
and taking expectation, we have
\begin{align*}
-\E\lc\dt\lp\sigma^2\Delta u\rp\sxns\rc
=-\E\lc\dt\lp\sigma^2\Delta u\rp\kxnk\rc+\tR_{6}(s),
\end{align*}
with $|\tR_{6}(s)|\leq Ch$ by Corollary \ref{cor}.
\subsubsection{The term $\tphi_7(s)=\bkn_s\partial\lp\sigma^2\Delta u\rp\sxns$}
Using \eqref{eq-betav} with $v=\partial(\sigma^2\Delta u)$,
integrating between $t_k$ and $s$,
taking expectation we have
\begin{align*}
-\E\lc\bkn_s\partial\lp\sigma^2\Delta u\rp\sxns\rc
=-\E\lc b\lp X^N\tk\rp\partial\lp\sigma^2\Delta u\rp\kxnk\rc+\tR_{7}(s),
\end{align*}
with $|\tR_{7}(s)|\leq Ch$ by Corollary \ref{cor}.
\subsubsection{The term	$\tphi_8(s)=\undeux\lv \gkn_s\rvd\Delta\lp\sigma^2\Delta u\rp\sxns$}
Using \eqref{eq-lvgammarvdv} with $v=\undeux\Delta(\sigma^2\Delta u)$,
integrating between $t_k$ and $s$,
and finally taking expectation we have
\begin{align*}
-\undeux\E\lc\lv \gkn_s\rvd\Delta\lp\sigma^2\Delta u\rp\sxns\rc
=-\undeux\E\lc\lv \sigma\lxnk\rvd\Delta\lp\sigma^2\Delta u\rp\kxnk\rc
+\tR_{8}(s),
\end{align*}
with $|\tR_{8}(s)|\leq Ch$ by Corollary \ref{cor}.
\subsection{Proof Theorem \ref{th-main} (i)}
The identity \eqref{eq-dec}
and the upper estimate in section 3.6 and 3.7 imply that
\begin{align}\E f\lp X^N_T\rp-\E f\lp X_T\rp
	=& \sum_{k=0}^{N-1}\E\intk\lb
	\int_t^{t\kk}\iNk(t\kk)ds
	+\undeux\int\tk^t\jNk(t_k)ds
	\rb dt +R\nonumber\\
	=&\undeux h^2\sum_{k=0}^{N-1}\E i^N_k(t\kk)
	+\unquart h^2\sum_{k=0}^{N-1}\E j^N_k(t_k)
	+R,\label{eq-EWR}\end{align}
where
$$R:=\sum_{k=0}^{N-1}\sum_{j=1}^7\int\tk^{t\kk}\int_t^{t\kk}R_j(s)dsdt
	+\sum_{k=0}^{N-1}\sum_{j=1}^{8}\int\tk^{t\kk}\int\tk^t\tR_j(s)dsdt.$$
Hence $\lv R\rv\leq Ch^2$.
Note that $\bkn\tkk=b\lp X^N\tkk\rp$.
Using \eqref{def-bzgn} and \eqref{eq-represent-z}
we deduce that $\zkn\tkk=\sigma\lp X^N\tk\rp\lp S_hb'\rp\lp X^N\tkk\rp$
and $\gkn\tkk=\sigma\lp X^N\tk\rp\lc1+h(S_hb')\lp X^N\tkk\rp\rc
=\sigma\lp X^N\tk\rp S_h\lp X^N\tkk\rp$.
Therefore, we deduce that
\begin{align*}
\iNk(t\kk)= &
\dt\lp b\partial u\rp\kkxnkk
	+\lc b\partial\lp b\partial u\rp\rc\kkxnkk
	+\undeux\sigma^2\lxnk \lc\Sh^2\Delta\lp b\partial u\rp\rc\kkxnkk
\\
&
-\lc b\dt\partial u\rc\kkxnkk
-\lc b^2\Delta u\rc\kkxnkk
-\sigma^2\lxnk \lc b'\Sh^2\Delta u\rc\kkxnkk
\\
&
-
\undeux\sigma^2\lxnk \lc b\Sh^2\partial^3u\rc\kkxnkk.
\end{align*}
Similary, $\gkn\tk=\sigma\lp X^N\tk\rp$.
So we have
\begin{align*}
\E j^N_k(t_k)
= &
\E\sigmad\dt\Delta u\kxnk
+
\E b\lxnkk\sigmad\partial^3u\kxnk
+
\undeux\E\sigma^4\partial^4u\kxnk
\\
&
+
2\E\Sh b'\lp X^N\tkk\rp\sigmad\Delta u\kxnk
-
\E\dt\lp \sigma^2\Delta u\rp\kxnk
\\
&
-
\E b\lxnkk\partial\lp\sigma^2\Delta u\rp\kxnk
-
\undeux\E\sigmad\Delta\lp\sigma^2\Delta u\rp\kxnk.
\end{align*}
Notice that $b$ and $\sigma$ do not depend upon t;
hence after simplification we have
\begin{align}
i^N_k(t\kk)
=&	b\partial\lp b\partial u\rp\kkxnkk
	-b^2\Delta u\kkxnkk
	+\undeux\sigmad\lxnk \lp \Sh^2b''\partial u\rp\kkxnkk\label{eq-iNk},
\end{align}
\begin{align}
\E\jNk(t_k)
=&	\E b\lxnkk\sigmad\partial^3u\kxnk
	+\undeux\E\sigma^4\partial^4u\kxnk
	-\E b\lxnkk\partial\lp\sigmad\Delta u\rp\kxnk\nonumber\\&
	-\undeux\E\sigmad\Delta\lp\sigmad\Delta u\rp\kxnk
	+2\E b'\Sh\lxnkk\sigmad\Delta u\kxnk\label{eq-jNk}.
\end{align}
Corollary \ref{cor} implies the existence of a constant $C$,
such that for all $k=0,\ldots,N-1$:
$$\lv\E\lp\iNk(t\kk)+\jNk(t_k)\rp\rv\leq C.$$
Using this bound with \eqref{eq-EWR} proves the first part of Theorem \ref{th-main}.
\subsection{Proof Theorem \ref{th-main} (ii)}
We at first prove the following lemma, 
which upper estimates the error in the approximation of
an integral by a Riemann sum.
\begin{lemma}
Let $v$ and $w$ in $C^{\infty,\infty}_b([0,T]\times\R)$.
Then there exists a constant $C$ independent of $h$
such that
$$\lv h\sum_{k=0}^{N-1}\E v\kkxnkk w\kxnk
	-\E\int_0^Tvw(t,X_t)dt\rv\leq Ch.$$
\end{lemma}
\proof Using \eqref{eq-vrxnr} 
multiply by $w(t_k,X^N\tk)$ and taking expected value,
we deduce
for $v\in C^{1,2}$,
$\E v\kkxnkk w\kxnk=\E vw\kxnk+A_k$, where
\begin{align*}
A_k
:=	\E w\kxnk\int\tk^{t\kk}\lb \dt v+\bkn_t\partial v
		+\undeux\lv\gkn_t\rvd\Delta v\rb \lp t,X^N_t\rp dt.
\end{align*}
This yields
$$h\sum_{k=0}^{N-1}\E v\kkxnkk w\kxnk
	-\E\int_0^Tvw(t,X_t)dt
=	\sum_{k=0}^{N-1}(hA_k+hB_k+C_k).$$
where
\begin{align*}
B_k
:=&	\E (vw)\kxnk-\E (vw)\lp t_k,X\tk\rp\\
C_k
:=&	h\E (vw)\lp t_k,X\tk\rp-\E\int\tk^{t\kk}(vw)(t,X_t)dt.
\end{align*}
Using the Cauchy-Schwarz inequality, the fact that
$\dt v$, $\partial v$ and $\Delta v$ have polynomial grow so that
Corollary \ref{cor} can be applied, we deduce
\begin{align*}
\lv A_k\rvd
\leq&	\E\lc\lv w\kxnk\rvd\rc\E\lc\lv \int\tk^{t\kk}
		\lb \dt v+\bkn_t\partial v+\undeux\lv\gkn_t\rvd\Delta v\rb
		\lp t,X^N_t\rp dt\rvd\rc\\
\leq&	Ch\E\int\tk^{t\kk}\lv
		\lb \dt v+\bkn_t\partial v+\undeux\lv\gkn_t\rvd\Delta v\rb
		\lp t,X^N_t\rp \rvd dt\\
\leq& 	Ch^2,
\end{align*}
Hence, $\lv A_k\rv\leq Ch$ which implies
$h\sum_{0\leq k\leq N-1}\lv A_k\rv\leq Ch$.\\
Since $(vw)(t_k,.)$ is in $C^{\infty}_b$, we use Theorem \ref{th-main} (i),
changing $T$ by $t_k$, which yields $\lv B_k\rv\leq Ch$
and then $h\sum_{0\leq k\leq N-1}\lv B_k\rv\leq Ch$.
Finally, It\^o's formula implies
\begin{align*}C_k=&\E\int\tk^{t\kk}\lb (vw)\lp t_k,X\tk\rp -(vw)\lp t,X_t\rp\rb dt\\
=&-\E\int\tk^{t\kk}\int\tk^t
\lb \dt(vw)+b\partial(vw)+\undeux\sigmad\Delta(vw)\rb(s,X_s)dsdt.
\end{align*}
Once more the polynomial growth imposed on $v$, $w$ and their partial derivatives implies that
$\lv C_k\rv\leq Ch^2$ and then
$\sum_{0\leq k\leq N-1}\lv C_k\rv\leq Ch$.
This concludes the proof.
\endproof
Now we introduce the function $\psih:[0,T]\times\R\rightarrow\R$ defined by
\begin{align}
\label{eq-psih}
\psih(t,x)
:=&	\undeux b\partial(b\partial u)(t,x)
	-\undeux b^2\Delta u(t,x)
	+\unquart\sigmad\Sh^2b''\partial u(t,x)
	+\unquart b\sigmad\partial^3u(t,x)\\&
	+{1\over 8}\sigma^4\partial^4u(t,x)
	+\undeux b'\Sh\sigmad\Delta u(t,x)
	-\unquart b\partial(\sigmad\Delta u)(t,x)
	-{1\over8}\sigmad\Delta(\sigmad\Delta u)(t,x)\nonumber
\end{align}
Using the expression of $\iNk$ in \eqref{eq-iNk} (resp. $\jNk$ in \eqref{eq-jNk})
and the previous lemma, we deduce
\begin{align}\label{eq-ijpsi}\lv\undeux h\sum_{k=0}^{N-1} \E\iNk(t\kk)
	+\unquart h\sum_{k=0}^{N_1}\E\jNk(t_k)
	-\int_0^T\E\psih(t,X_t)dt\rv
\leq Ch.
\end{align}
Using the definitions of $\psi_i$ and $\psih$
given in \eqref{eq-psii} and \eqref{eq-psih} respectively, we have
\begin{align*}
\psih(t,x)-\psi_i(t,x)
=&	\unquart\lb\sigmad\Sh^2b''\partial u
	+b\sigmad\partial^3u+2b'\Sh\sigmad\Delta u
	-\sigmad\Delta(b\partial u)\rb(t,x)\\
=&	\unquart\sigmad\lp\Sh^2-1\rp b''\partial u(t,x)
	+\undeux b'\lp \Sh-1\rp\sigmad\Delta u(t,x).
\end{align*}
Since $\lp\Sh-1\rp(x)={hb'\over 1-hb'}(x)$
and $|\Sh(x)+1|\leq C$ for $h\in(0,h^*)$, we have
$\lv\lp\Sh-1\rp(x)\rv+\lv\lp\Sh^2-1\rp(x)\rv\leq Ch$,
where as usually $C$ does not depend on $N$ and $h$.
This yields
$$\lv\int_0^T\E\lb\psih(t,X_t)-\psi_i(t,X_t)\rb dt\rv\leq Ch.$$
This last equation with \eqref{eq-EWR} and \eqref{eq-ijpsi} concludes the proof.

\subsection*{Acknowledgments:} 
The author wishes to thank Annie Millet for
many helpful comments.

\end{document}